\documentclass[11pt,reqno]{amsart}
\usepackage[utf8x]{inputenc}
\usepackage{mathpazo} 
\usepackage{ucs}
\usepackage[english]{babel}
\usepackage{amsmath,slashed}
\usepackage{amsfonts}
\usepackage{mathtools}
\usepackage{amssymb}
\usepackage{subcaption}
\usepackage[shortlabels]{enumitem}
\usepackage{bbm}
\usepackage{graphicx}
\usepackage[left=2cm,right=2cm,top=2cm,bottom=2cm]{geometry}

\newlength{\intwidth}

\usepackage{amssymb}
\usepackage[centertags]{amsmath}
\usepackage{verbatim}

\newcommand\om\omega

\newcommand{\si}{\sigma}

\def\RR{\mathbb{R}}

\newcommand{\foralll}{\forall\hspace{1mm}}

\renewcommand\ell{{ l }}
\DeclareMathOperator\Real{Re}

\def\TT{\mathbb{T}}

\newcommand\vp\varphi
\newcommand\ka\kappa
\newcommand\te\theta
\newcommand\vka\varkappa

\newcommand\de\delta
\newcommand\La\Lambda
\newcommand\la\lambda
\newcommand\ga\gamma

\newcommand\Ga\Gamma

\newcommand\cJ{\mathcal J}

\newcommand{\triple}[1]{{\left\vert\kern-0.25ex\left\vert\kern-0.25ex\left\vert #1
        \right\vert\kern-0.25ex\right\vert\kern-0.25ex\right\vert}}

\usepackage[toc, page]{appendix}

\newcommand{\norm}[1]{\left\lVert#1\right\rVert}

\newcommand{\bR}{\mathbb{R}}
\newcommand{\bN}{\mathbb{N}}
\newcommand{\bE}{\mathbb{E}}

\newcommand{\bZ}{\mathbb{Z}}

\makeatletter
\renewcommand*{\@fnsymbol}[1]{\ensuremath{\ifcase#1\or *\or \star\or ***\or
   \mathsection\or \mathparagraph\or \|\or **\or \dagger\dagger
   \or \ddagger\ddagger \else\@ctrerr\fi}}
\makeatother

\usepackage{fancyhdr}
\newcommand\ep\varepsilon

\newtheorem{theorem}{Theorem}[section]
\newtheorem{lemma}[theorem]{Lemma}

\newtheorem{proposition}[theorem]{Proposition}
\newtheorem{corollary}[theorem]{Corollary}

\theoremstyle{definition}
\newtheorem{remark}[theorem]{Remark}

\numberwithin{equation}{section}

\renewcommand\leq\leqslant
\renewcommand\geq\geqslant

\usepackage[pdftex]{hyperref}
\hypersetup{
	colorlinks=true,
	linkcolor=blue,
	citecolor=blue,
	filecolor=blue,
	urlcolor=blue,
}

\title[Integral representations and asymptotic expansions for Bessel series]{Integral representations and asymptotic expansions for the second type Neumann series of Bessel functions of the first kind}

\author{\'Alvaro Romaniega}
\address{Instituto de Ciencias Matem\'aticas, Consejo Superior de
 Investigaciones Cient\'\i ficas, 28049 Madrid, Spain}
\email{alvaroromaniega@gmail.com}

\begin{document}
\maketitle

\begin{abstract} In this paper we study the following Bessel series $\sum _{l=1}^{\infty } {J_{l+m'}(r)J_{l+m}(r)}{(l+\beta)^\alpha}$ for any $m,m'\in\bZ$, $\alpha\in\bR$ and $\beta>-1$. They are a particular case of the second type Neumann series of Bessel functions of the first kind. More specifically, we derive fully explicit integral representations and study the asymptotic behavior with explicit terms. As a corollary, the asymptotic behavior of series of the derivatives of Bessel functions can be understood. 
\end{abstract}

\section{Introduction.}
Bessel functions of the first kind and order $\nu$ can be defined as
$$
J_{\nu}(x)=\sum_{n=0}^\infty\frac{(-1)^n}{n!\Gamma(n+\nu+1)}\left(\frac{x}{2}\right)^{2n+\nu}\,.
$$
As is well-known, they are a  particular solution of the Bessel differential equation,
\begin{equation*}
	x^2y''(x)+xy'(x)+(x^2-\nu^2)y(x)=0\,.
\end{equation*}
Bessel functions are ubiquitous functions in (applied) mathematics and theoretical physics. In particular, these series of Bessel functions appear frequently and they have been studied in detail, \cite[Chapters XVI-XIX]{Wat95}.

In this paper we are interested in the second type Neumann series of Bessel functions of the first kind, \cite[Section 2.5]{BMP17},
\begin{equation}\label{eq:Neumann series}
	\mathfrak{N}^{a,b}_{\mu,\nu}(x)\coloneqq\sum_{n\ge 1}\alpha_n J_{\mu+an}(x)J_{\nu+bn}(x)\,,	
\end{equation}
being all parameters real numbers. Classical examples of these series are von Lommel’s series and the Al-Salam series, see \cite[Section 2.5]{BMP17} and references therein. Here we want to study the case of $a=b=1$ and
$$
\alpha_n = (n+\beta)^{\alpha}\,,
$$
for $\alpha\in\bR$ and $\beta>-1$. That is\footnote{We have used $m,m'$ instead of the $\mu,\nu$ of (\ref{eq:Neumann series}). Hereafter, $\mu\coloneqq m+m'$, $\nu\coloneqq m-m'$, see below.},
\begin{equation}\label{eq:def S}
	S_{\alpha,\beta,m,m'}(r)\coloneqq\sum _{l=1}^{\infty } {J_{l+m'}(r)J_{l+m}(r)}{(l+\beta)^\alpha}.
\end{equation}
Classical examples of these series are consequences of Neumann’s Addition Theorem, see \cite[(10.23.3),(10.23.4)]{Olv10},
\begin{align*}
	\sum _{k=0}^{2 n} (-1)^k J_k(r) J_{2 n-k}(r)+2 \sum _{l=1}^{\infty } J_l(r) J_{l+2 n}(r)&=0,\\ {J_{0}}^{2}\left(z\right)+2\sum_{k=1}^{\infty}{J_{k}}^{2}\left(z\right)&=1,
\end{align*}
for $n\ge1$. Furthermore, from this theorem one can easily derive exact formulas for derivatives of Bessel functions, see \cite[Appendix B]{EPR21}. For instance,
\begin{align*}
	\sum_{l=0}^\infty{\epsilon_l} J_l'(r)^2&=\frac12\,,\quad
	\sum_{l=0}^\infty{\epsilon_l} l^2J_l(r)J_l'(r)=\frac r4\,,\quad
	\sum_{l=0}^\infty{\epsilon_l} l^4J_l(r)^2=\frac{r^2(4+3r^2)}8\,,
\end{align*}
where $\epsilon_l:=2-\de_{l,0}$ is Neumann's factor. Another classical example is the following Tur\'an type inequality for all
$x\in\mathbb{R}$ and $\nu>-1$, \cite[p. 384]{TN51},
\begin{equation}\label{turan1}
	\Delta_{\nu}(x)= J_{\nu}^2(x)-J_{\nu-1}(x)J_{\nu+1}(x)\geq0\,,
\end{equation}
where\footnote{It might seem that this is not of the type \eqref{eq:def S}, but it can be easily converted to a linear combination using, for instance, $$\frac{1}{(\beta +n) (\beta +n+2)}=\frac{1}{2 (\beta +n)}-\frac{1}{2 (\beta +n+2)}.$$}
\begin{equation} \label{turan4}
	\Delta_{\nu}(x)\coloneqq\frac{1}{\nu+1}J_{\nu}^2(x)+\frac{2}{\nu+2}J_{\nu+1}^2(x)
	+2\nu\sum_{n\geq2}\frac{J_{\nu+n}^2(x)}{(\nu+n-1)(\nu+n+1)}.
\end{equation}

Recently, these series have appeared in the study of critical points of random fields satisfying the Helmholtz equation on the plane, \cite{EPR21}. This is expected as, by definition, 
the random field~$u$ is
\begin{equation}\label{defu2}
	u := \sum_{l} a_l \,\si_l\, e^{i l\te}\, J_l(r)\,,\qquad \si_l:= \begin{cases}|l|^{-s} & \text{if }l\neq0\,,\\ 0& \text{if }l=0\,,\end{cases}
\end{equation}
where the real and imaginary parts of $a_l$
are independent standard Gaussian random variables subject to the
constraint $a_l=(-1)^l\overline{a_{-l}}$ (which makes~$u$ real valued),
$(r,\te)\in \RR^+\times\TT$ are the polar coordinates. Therefore, the (co)variance kernel of the random function~\eqref{defu2} is
\[
K(r,\te):=\bE[u(r,\te)\,u(r,\te)]=4\sum _{l=1}^\infty l^{-2s} J_l(r) J_l(r)\,,
\]
see \cite[Remark 4.2]{EPR21}. The covariance kernel is one of the main inputs in the Kac-Rice formula, e.g. \cite{AW09}, which can be used to compute the expected value of the number of critical points in a given region. There, a method to compute the asymptotics of $\cJ_{s,m,m'}(r)\coloneqq S_{-2\alpha,\beta=0,m,m'}(r)$ is given. The basic idea of the strategy is to decompose the series as follows

$$\cJ_{s,m,m'}(r)=\cJ^{\textnormal{I}}_{s,m,m'}(r)+\cJ^{\textnormal{II}}_{s,m,m'}(r)+\cJ^{\textnormal{III}}_{s,m,m'}(r)\,,$$
where (for some $\delta>0$ small):
\begin{itemize}
	\item $\textnormal{I}$ only involves ``frequencies'' $\lambda\coloneqq l/r$ smaller than  $(1-\de)$, 
	\item $\textnormal{II}$
	involves $\lambda$ close to~1 (more precisely, $|\lambda-1|<2\delta$), and
	\item $\textnormal{III}$ involves $\lambda$ larger than $(1+\de)$.	
\end{itemize}\vspace{-2mm}
Then, different oscillatory techniques are used in each region to extract the leading term and control the errors. This is done for every $\alpha$, using different techniques for different regions of this parameter.

A completely different approach is followed here. First, we divide the series between $\alpha\ge0$ and $\alpha<0$ and derive integral representations in the latter case. Once we have those integral representations, we use different asymptotic techniques from harmonic analysis (such as the stationary phase method and the Hankel transform) to compute their asymptotic expansions. Finally, the case of $\alpha\ge0$ is treated using a ``trick'' (Lemma \ref{lemma:trick}) to reduce it to the previous case. This method allow us to obtain both integral representations and asymptotic expansions.

Beyond the fact that the methods used are different, we believe there are also some worth noting improvements. First, the series considered are more general as the parameter $\beta$ can take any value in $(-1, +\infty)$, not just $0$. Second, we also provide integral representations, which is an important result by itself and can be useful for other purposes apart from asymptotic expansions when $r\to\infty$. See \cite[Section 2.5]{BMP17} and references therein for some literature on other integral representations for these types of series. Third, our method allows us to compute the whole asymptotic expansion, not just the leading term. For the sake of simplicity, we have basically limited our attention to the first terms, but using Remark \ref{rem:higher order terms} and \ref{rem:higher order terms integer} obtaining higher order terms is straightforward. In fact, in Proposition \ref{prop:asymp exp alpha integer} second order terms are easily computed for convenience (they are needed for the case of $\alpha=0$). This also results in a better understanding of the asymptotics for some series of derivatives of Bessel functions. Comparing \cite[Corollary 3.7]{EPR21} with Corollary \ref{cor:deriv series alpha -1} here, we can see, e.g., the following improvement,
\begin{align*}
		\sum _{l=1}^{\infty } \frac{J_l(r) J'_l(r)}{l}= O\left(r^{-1}\right)\Rightarrow\sum _{l=1}^{\infty } \frac{J_l(r) J'_l(r)}{l+\beta}= -\frac{\Phi (-1,1,\beta +1) \cos (2 r)}{\pi  r}+o(r^{-1})\,,
\end{align*}
where $\Phi$ is the Lerch transcendent function.
The organization and main results of this paper are as follows. Relevant notation is introduced in Section \ref{sec:notation}. In Section  \ref{sec:int rep neg alpha} we derive the integral representations for $\alpha<0$. In particular, for $\mu\in\mathbb{Z}_{\ge0}$, we define
\begin{equation*}
	F_{\alpha,\beta,\mu}(\varphi)\coloneqq-\frac{1}{2} e^{-i \varphi  (\mu +2)} \left(e^{2 i \varphi  (\mu +2)} \Phi \left(-e^{2 i \varphi },\alpha ,\beta +1\right)+\Phi \left(-e^{-2 i \varphi },\alpha ,\beta +1\right)\right)\,,
\end{equation*}
where, as above, $\Phi$ is the Lerch transcendent function defined in Section \ref{sec:notation}. 
\begin{theorem}\label{prop:integral rep} Let $\alpha,\beta+1>0$ and $m,m'=0,1,\ldots$, then the series \eqref{BesselSeries} can be expressed in the following integral forms:
	\begin{equation}\label{eq:integral rep exp}
		\sum _{l=1}^{\infty } \frac{J_{l+m'}(r)J_{l+m}(r)}{(l+\beta)^\alpha}=\frac{2i^{-\mu}}{\pi^2}\int_{0}^{\pi/2}\int_{0}^{\pi}e^{i2r\cos\varphi\cos\theta}F_{\alpha,\beta,\mu}(\varphi)\cos(\nu\theta) d\varphi d\theta\,,
	\end{equation}
	as a \textnormal{{two-dimensional exponential oscillatory integral}} and
	\begin{equation}\label{eq:hankel}
		\sum _{l=1}^{\infty } \frac{J_{l+m'}(r)J_{l+m}(r)}{(l+\beta)^\alpha}=\frac{(-1)^{m'}}{\pi} \int_{0}^{\pi}  \, J_{\nu}(2r\cos\varphi) F_{\alpha,\beta,\mu}(\varphi) d\varphi\,,
	\end{equation}	
	as a \textnormal{{one-dimensional Hankel transform}}, for $r\in\mathbb{R}$ where $\mu\coloneqq m+m'$ and $\nu\coloneqq m-m'$. 
\end{theorem} 
\noindent Using both integral representations, in Section \ref{sec: asymptotics neg alpha} we are able to compute its asymptotic expansion. In particular, 
\begin{theorem} Let $\mu\coloneqq m+m'$ and $\nu\coloneqq m-m'$. If $\alpha>0$ and $\alpha\notin \bN$, then there are constants $c$ such that
	\begin{equation*}
		\sum _{l=1}^{\infty } \frac{J_{l+m'}(z)J_{l+m}(z)}{(l+\beta)^\alpha}=\frac{c^1_{\alpha,\nu}}{r^\alpha}+\frac{c^2_{\alpha, \beta,\nu}+c^3_{\alpha,\beta}\sin \left(2r-{\pi\mu}/{2}\right)}{r}+O\left(\frac{1}{r^{\gamma}}\right)\,,
	\end{equation*}
with $\gamma\coloneqq \min\{\alpha+1, 2\}$, see Corollary \textnormal{\ref{cor:series exp alpha pos}} for explicit constants. Now let $\alpha\in\bN$. Similarly for $\alpha=1$,
	\begin{align*}
	\sum _{l=1}^{\infty } \frac{J_{l+m'}(z)J_{l+m}(z)}{(l+\beta)}=\frac{c^4_\nu \log r}{r}+
	\frac{c^5_{\beta,\nu}+c^3_{1,\beta}\sin \left(2r-{\pi\mu}/{2}\right)}{r}	+o\left(r^{-1}\right)\,,
	\end{align*}
and for $\alpha>1$,
\begin{align*}
	\sum _{l=1}^{\infty } \frac{J_{l+m'}(z)J_{l+m}(z)}{(l+\beta)^\alpha}=&\frac{c^2_{\alpha, \beta,\nu}+c^3_{\alpha,\beta}\sin \left(2r-{\pi\mu}/{2}\right)}{r}+O\left(\frac{1}{r^{2-\gamma}}\right)\,,
\end{align*} 
with $\gamma>0$ arbitrarily small, see Corollary \textnormal{\ref{cor:series expansion alpha integer}} for explicit constants.
\end{theorem}
\noindent In Section \ref{sec: int rep asymp pos alpha}, first we derive the asymptotics for $S_{\alpha,\beta,\mu,\nu}(r)$ when $\alpha\ge0$. In particular,
\begin{theorem}\label{prop:series expansion S} For $\alpha\ge0$ and $\beta,\mu,\nu$ as above,
	\begin{equation*}
		\sum _{l=1}^{\infty } {J_{l+m'}(z)J_{l+m}(z)}{(l+\beta)^\alpha}={c^1_{-\alpha,\nu}}r^\alpha+o(r^\alpha)=\frac{2^{-\alpha -1} \Gamma (\alpha +1) r^{\alpha }}{\Gamma \left(\frac{1}{2} (-\nu+\alpha +2)\right) \Gamma \left(\frac{1}{2} (\nu+\alpha +2)\right)}+o(r^\alpha)\,,
	\end{equation*}
where the coefficient must be understood as the continuous extension for the values of $\alpha$ where it is undefined.
\end{theorem}
In Proposition \ref{prop:int rep alpha nonneg} integral representations for this case are given. Finally, in Section \ref{sec: deriv series}, we derive asymptotic formulas for series of derivatives of Bessel functions.
\section{Notation.}\label{sec:notation}
With $\Phi$ we denote the Lerch transcendent function, i.e., 
\begin{equation*}
	\Phi (z,\alpha ,\beta )=\sum _{n=0}^{\infty } \frac{z^n}{(\beta +n)^{\alpha }}\,,
\end{equation*}
for $|z|<1$ (or $|z|\le1$ if $\alpha>1$), see \cite[25.14.1]{Olv10}, and for other values it is extended by analytic continuation. As a particular case we have $\text{Li}_{\alpha }$, the polylogarithm of order $\alpha$. The function $\text{Li}_{\alpha}(z)$ is defined as a series
\begin{equation}\label{eq:def Li}
	\mathrm{Li}_{\alpha}\left(z\right)=\sum_{n=1}^{\infty}\frac{z^{n}}{n^{\alpha}}\,,
\end{equation}
for $|z|<1$ (or $|z|\le1$ if $\alpha>1$), see \cite[25.12.10]{Olv10}, and for other values it is extended by analytic continuation. We have that
\begin{equation*}
	\mathrm{Li}_{\alpha}\left(z\right)=z\Phi\left(z,\alpha,1\right)\,.
\end{equation*}
Also, the family of Struve functions is denoted by $\textnormal{\pmb{H}}_{\nu }(z)$. Recall that Struve functions are particular solutions of the non-homogeneous Bessel's differential equation
$$
\frac{{ d}^{2}y}{{ dz}^{2}}+\frac{1}{z}\frac{ dy}{  dz}+\left(1-\frac{\nu^{2}}{z^{2}}\right)y=\frac{(\tfrac{1}{2}z)^{\nu- 1}}{\sqrt{\pi}\Gamma\left(\nu+\tfrac{1}{2}\right)}
$$
defined by the power-series expansion
$$
\textnormal{\pmb{H}}_{\nu }(z)\coloneqq(\tfrac{1}{2}z)^{\nu+1}\sum_{n=0}^{\infty}\frac{(-1)^{n}(\tfrac{ 1}{2}z)^{2n}}{\Gamma\left(n+\tfrac{3}{2}\right)\Gamma\left(n+\nu+\tfrac{3}{2} \right)}.
$$
The $(p,q)$-regularized hypergeometric function will be denoted by $_p\tilde{F}_q$.

\section{Integral representations for Bessel series}\label{sec:int rep neg alpha}
\begin{figure}[h]
\begin{subfigure}{0.495\textwidth}
		\centering
	\includegraphics[width=\linewidth]{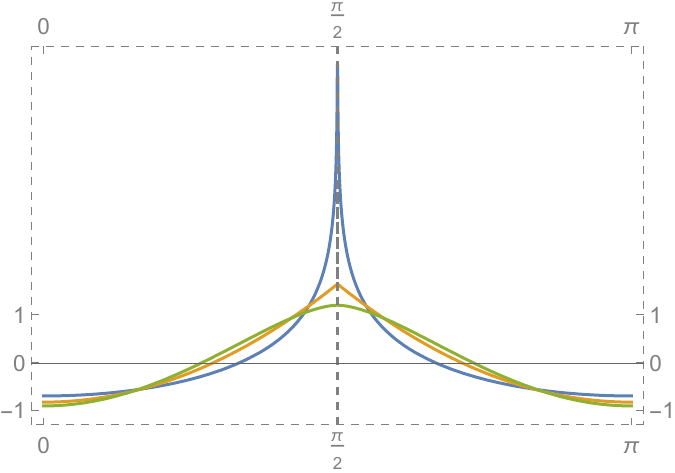}
	\caption{Representation of $F_{\alpha,0}$ for several values of $\alpha$.}
	\label{fig:fig1}
\end{subfigure}
	\begin{subfigure}{0.495\textwidth}
		\centering
	\includegraphics[width=\linewidth]{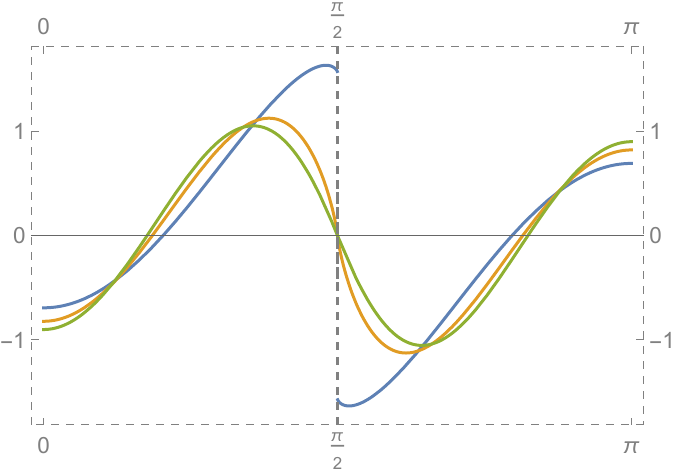}
	\caption{Representation of $F_{\alpha,1}$ for several values of $\alpha$. }
	\label{fig:fig2}
	\end{subfigure}
\begin{subfigure}{0.495\textwidth}
	\centering
	\includegraphics[width=\linewidth]{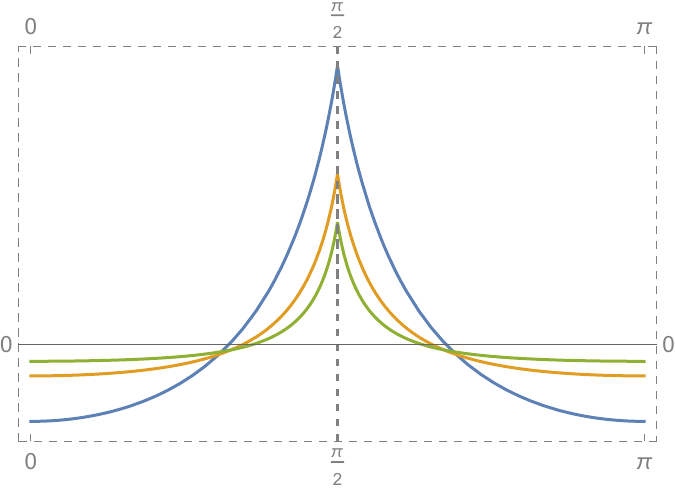}
	\caption{Representation of $F_{2,\beta,0}$ for several values of $\beta$.}
	\label{fig:fig3}
\end{subfigure}
\begin{subfigure}{0.495\textwidth}
	\centering
	\includegraphics[width=\linewidth]{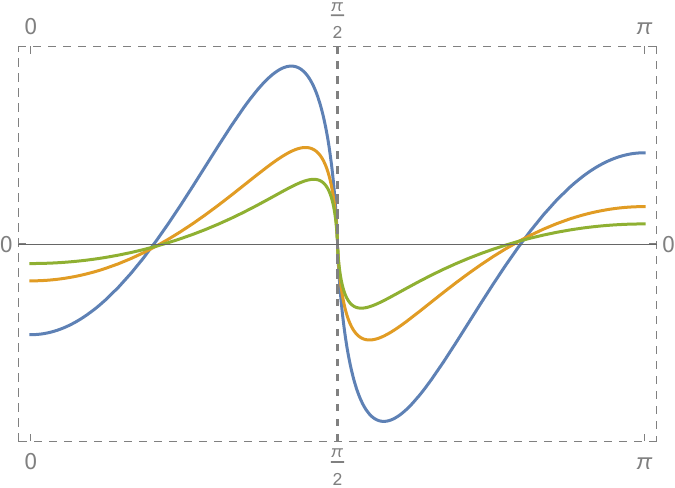}
	\caption{Representation of $F_{2,\beta,0}$ for several values of $\beta$.}
	\label{fig:fig4}
\end{subfigure}
\caption{Representations of $F_{\alpha,\beta,\mu}$ for different parameters. The dashed line represents $\varphi=\pi/2$. }
\label{fig:rep F}
\end{figure}

The purpose of this section is give an integral representation of
\begin{equation}\label{BesselSeries}
S_{-\alpha,\beta,m,m'}(r)\coloneqq\sum _{l=1}^{\infty } \frac{J_{l+m'}(r)J_{l+m}(r)}{(l+\beta)^\alpha}\,,
\end{equation}
where $m,m'=0,1,\dots$, $\alpha>0$ and $\beta>-1$ so that $l+\beta>0$. For this we need some definitions. We define, for $\mu\in\mathbb{Z}_{\ge0}$,
\begin{equation*}
	F_{\alpha,\beta,\mu}(\varphi)\coloneqq-\frac{1}{2} e^{-i \varphi  (\mu +2)} \left(e^{2 i \varphi  (\mu +2)} \Phi \left(-e^{2 i \varphi },\alpha ,\beta +1\right)+\Phi \left(-e^{-2 i \varphi },\alpha ,\beta +1\right)\right)\,,
\end{equation*}
where $\Phi$ is the Lerch transcendent function defined in Section \ref{sec:notation}. A particular interesting situation is
\begin{equation*}
F_{\alpha,m}(\varphi)\coloneqq F_{\alpha,0,m}(\varphi)= \frac{1}{2} \left(e^{ i \varphi  m} \text{Li}_{\alpha }\left(-e^{2 i \varphi }\right)+ e^{-i \varphi  m}\text{Li}_{\alpha }\left(-e^{-2 i \varphi }\right)\right)\,,
\end{equation*}
where $\text{Li}_{\alpha }$ is the polylogarithm of order $\alpha$ defined in Section \ref{sec:notation}. Before the main proposition of this section, let us study the function $F_{\alpha,m}$, in particular, its series expansion.
\begin{lemma}\label{lemma:F alpha m} For $0\le\varphi<\pi/2$, $\alpha>0$ and $\mu\in\mathbb{Z}_{\ge0}$, we have
	\begin{equation}
		F_{\alpha,\beta,\mu}(\varphi)=\sum_{l=1}^{\infty } \frac{(-1)^{l} \cos (\varphi (\mu+2 l))}{(l+\beta)^{\alpha }}=(-1)^\mu F_{\alpha,\beta,\mu}(\pi-\varphi)\,,
	\end{equation}
where the equality is true for $\varphi=\pi/2$ if $\alpha>1$. $F_{\alpha,\beta,\mu}$ is a smooth function for $0\le\varphi<\pi/2$ and, possibly, with singularities at $\varphi=\pi/2$ if $\alpha\le 1$. The same singularities can appear for the derivatives if $\alpha>1$.
\end{lemma}
\begin{proof}
	By \eqref{eq:def Li} and the comments below, the only non-trivial case to consider is $0<\alpha\le 1$ as the values of our interest ($z=-e^{i2\varphi}$) are on the boundary of the disk of convergence of the series, i.e., $z\in \partial D\backslash D$ where $D$ is the disk of convergence. First note that the following geometric sum gives
	\begin{align}\label{eq:F alpha=0}
	A_N(\rho,\varphi)\coloneqq	\sum _{n=1}^N (-1)^n \rho^n e^{i \varphi2n}=\frac{\rho  e^{2 i \varphi } \left(-1+(-\rho )^N e^{2 i N \varphi }\right)}{1+\rho  e^{2 i \varphi }}\,,
	\end{align}
	where $\rho\in[0,1]$. Thus  
	\begin{equation*}
		|A_N(\rho,\varphi)|\le \frac{2}{|1+\rho e^{2 i \varphi }|}\le M_\delta
	\end{equation*}	
	for $\varphi\in[-\frac{\pi}{2}+\delta,\frac{\pi}{2}-\delta], ~ \foralll\delta>0$. Obviously this constant blows up as $\delta\to0$. Hence, if
	\begin{equation*}
 	S_N(\rho,\varphi)\coloneqq \sum _{n=1}^N\frac{ i^{-2 n} \rho^n e^{i \varphi 2n}}{(n+\beta)^\alpha}\,,
	\end{equation*}  
	by summation by parts for $n>m$, 
	\begin{align*}
		|S_n-S_m|(\rho,\varphi)&=\left|\frac{A_n(\rho,\varphi)}{(n+\beta)^\alpha}-\frac{A_m(\rho,\varphi)}{(m+\beta)^\alpha}+\sum_{k=m}^{n-1}A_k(\rho,\varphi)\left(\frac{1}{(k+\beta)^\alpha}-\frac{1}{(k+1+\beta)^\alpha}\right)\right|\le\\
		&\le M_\delta\left(\frac{1}{(n+\beta)^\alpha}+\frac{1}{(m+\beta)^\alpha}+\frac{1}{(m+\beta)^\alpha}-\frac{1}{(n+\beta)^\alpha}\right)=\frac{2M_\delta}{(m+\beta)^\alpha},
	\end{align*}
	where the last inequality follows from a telescopic cancellation. As the RHS goes to zero as $m\to\infty$, the sequence is Cauchy and
	\begin{equation*}
		|S(\rho,\varphi)-S_m((\rho,\varphi))|< \varepsilon
	\end{equation*}
	for $m>\left(2M_\delta/\varepsilon\right)^\frac{1}{\alpha}-\beta$, independent of $\rho,\varphi$. Thus, the convergence is uniform so $S(\rho,\varphi)$ is continuous if $\varphi\neq \frac{\pi}{2}$. For $\rho<1$,  $S(\rho,\varphi)=\rho  \left(-e^{2 i \varphi }\right) \Phi \left(-e^{2 i \varphi } \rho ,\alpha ,\beta +1\right)$ (by the definition of $\Phi$ as a series). Also, for $\varphi\in[-\frac{\pi}{2}+\delta,\frac{\pi}{2}-\delta]$ we have $S(1,\varphi)=  -e^{2 i \varphi } \Phi \left(-e^{2 i \varphi } ,\alpha ,\beta +1\right)$ (understood as the analytic continuation of $\Phi$ for $\rho<1$). As $\delta>0$ is arbitrary, the proof of the series expansion follows straightforwardly taking the real part, i.e.,
		\begin{align*}
			\sum _{l=1}^{\infty } \frac{(-1)^{l} \cos (\varphi (\mu+2 l))}{(l+\beta)^{\alpha }}&=\sum _{l=1}^{\infty } \frac{\Real\left(e^{i\varphi\mu}(-e^{i2\varphi})^l\right)}{(l+\beta)^{\alpha }}=\\
			&=\Real\left(-e^{i \varphi (\mu +2)} \Phi \left(-e^{2 i \varphi },\alpha ,\beta +1\right)\right)=F_{\alpha,\beta,m}(\varphi)\,.
		\end{align*}
	The parity property follows straightforwardly from the series expansion.
	
	We have the following integral representation \cite[25.14.5]{Olv10}
	\begin{equation}\label{IntRepLi}
	\Phi \left(-e^{2 i \varphi },\alpha ,\beta +1\right)=\frac{1}{\Gamma (\alpha )}\int_0^{\infty } \frac{t^{\alpha -1} e^{-t(\beta +1)}}{1+e^{-t+2 i \varphi }} \, dt\,,
	\end{equation}
	so we might have singularities at $\varphi=\frac{\pi}{2}$ if $\alpha\le1$ and the function is $C^\infty$ except at $\varphi=\frac{\pi}{2}$. This is even true if $\alpha>1$ as
	$$
	\frac{\partial \text{Li}_{\alpha }\left(-e^{\pm 2 i \varphi }\right)}{\partial \varphi}=\pm 2 i \text{Li}_{\alpha -1}\left(-e^{\pm 2 i \varphi }\right),
	$$
	\begin{equation*}
		\frac{\partial \Phi \left(-e^{2 i \varphi },\alpha ,\beta +1\right)}{\partial \varphi}=2 i \left(\Phi \left(-e^{2 i \varphi },\alpha -1,\beta +1\right)-(\beta +1) \Phi \left(-e^{2 i \varphi },\alpha ,\beta +1\right)\right)\,.
	\end{equation*}
\end{proof}
See Figure \ref{fig:rep F} for some representations.
\begin{remark}\label{rem:convergence Li alpha=0}
	Obviously the analytic continuation and the sum can differ on the boundary (the region of interest), that is why we need to check that they coincide. For instance, taking $\beta=0$ for simplicity, this is clear for $\alpha=0$ where the partial sums are given by \eqref{eq:F alpha=0} and they are not convergent for $\rho=1$. Nevertheless, the analytic continuation is well-defined and equals
	$$
	\mathrm{Li}_{0}(-e^{2i\varphi})=\frac{-e^{i (\varphi  (m+2))}}{1+ e^{2 i \varphi }}.
	$$
\end{remark}
We are ready to prove the main theorem of this section.
\begin{proof}[Proof of Theorem \ref{prop:integral rep}]
First, we know that for $z\in\mathbb{C}$ \cite[10.9.26]{Olv10}
\begin{equation*}
J_{\mu}\left(z\right)J_{\nu}\left(z\right)=\frac{2}{\pi}\int_{0}^{\pi/2}J_{\mu%
	+\nu}\left(2z\cos\theta\right)\cos\left((\mu-\nu)\theta\right) d\theta\,,
\end{equation*}
if $\Real(\mu+\nu)>-1.$
Second, we also know that \cite[10.9.2]{Olv10}
\begin{equation}\label{eq:integral rep Jn}
J_{n}\left(z\right)=\frac{1}{\pi}\int_{0}^{\pi}\cos\left(z\sin\varphi-n\varphi
\right) d\varphi=\frac{i^{-n}}{\pi}\int_{0}^{\pi}e^{iz\cos\varphi}\cos
\left(n\varphi\right) d\varphi\,,
\end{equation}
for $n\in\mathbb{Z}$. From these two expressions we obtain
$$
J_{l+m'}\left(z\right)J_{l+m}\left(z\right)=\frac{2i^{-\mu}}{\pi^2}\int_{0}^{\pi/2}\int_{0}^{\pi}e^{i2z\cos\varphi\cos\theta}(-1)^l\cos
\left((2l+\mu)\varphi\right)\cos\left(\nu\theta\right) d\varphi d\theta\,.
$$
Thus, if we can interchange the integral and the series we will have:
\begin{equation*}
\sum _{l=1}^{\infty } \frac{J_{l+m'}(z)J_{l+m}(z)}{(l+\beta)^\alpha}=\frac{2i^{-\mu}}{\pi^2}\int_{0}^{\pi/2}\int_{0}^{\pi}e^{i2z\cos\varphi\cos\theta}F_{\alpha,\beta,\mu}(\varphi)\cos(\nu\theta) d\varphi d\theta,
\end{equation*}
where we have used Lemma \ref{lemma:F alpha m}. To justify the swap, note that $F_{\alpha,\beta,\mu}$ can be written as
\begin{equation*}
F_{\alpha,\beta,\mu}(\varphi)=\frac{1}{2}e^{-i \varphi \mu}\sum _{l=1}^{\infty } (-1)^l (l+\beta)^{-\alpha } e^{-i 2l\varphi }+\frac{1}{2}e^{+i \varphi \mu} \sum _{l=1}^{\infty }(-1)^l (l+\beta)^{-\alpha } e^{i 2l \varphi},
\end{equation*}
i.e., as a Fourier series. Hence, the swap will be justified if we can ensure this series converges in $L^p$ for some $p\ge 1$. Indeed, if $S_N$ denotes the $N$-partial sum of the Fourier series, 
\begin{align}\label{eq:swap int sum}
&\left|\int_{0}^{\pi/2}\int_{0}^{\pi}e^{i2r\cos\varphi\cos\theta}\left(\sum _{l=N+1}^{\infty }\frac{(-1)^l\cos
\left((2l+\mu)\varphi\right)}{(l+\beta)^\alpha}\right)\cos(\nu\theta) d\varphi d\theta\right|\lesssim\\ &\hspace{7cm}\lesssim\norm{F_{\alpha,\beta,\mu}-S_N}_{L^1([0,\pi])}\lesssim\norm{F_{\alpha,\beta,\mu}-S_N}_{L^p([0,\pi])}\,,\nonumber
\end{align}
for $p>1$. This is clear\footnote{It is also obvious by Fubini if $\alpha>1$.} for $p=2$ and $\alpha>1/2$ as then
\[
\sum _{l=1}^{\infty }\frac{1}{(l+\beta)^{2\alpha}}<\infty\,,
\]
but a more detailed analysis of the function, not only the coefficients, is needed for $\alpha\in(0,1/2)$ (the function might have singularities at $\varphi=\pi/2$ by Lemma \ref{lemma:F alpha m}). From the series representation, \cite[p. 29]{Bat53}, 
\begin{align}\label{eq:series LerchT}
\Phi (z,s,a)=z^{-a} \Gamma (1-s) \log ^{s-1}\left(\frac{1}{z}\right)+z^{-a} \left(\sum _{k=0}^{\infty } \frac{\log ^k(z) \zeta (s-k,a)}{k!}\right) \text{ if }\neg (s\in \mathbb{Z}\land s>0)
\end{align}
for $|\log z|<2\pi$ where $\zeta$ is the (analytic continuation) of the Hurwitz zeta function, we can conclude (after a lengthy but straightforward computation) 
\begin{equation*}
	|F_{\alpha,\beta,\mu}(\varphi)|\sim_{\alpha,\mu} |\pi/2-\varphi|^{\alpha-1}\,.
\end{equation*}
More precisely,
\begin{equation*}
F_{\alpha,\beta,\mu}(\varphi)\sim\begin{cases}
\frac{1}{2} \Gamma (1-\alpha ) \left(-(\pi -2 \varphi )^{\alpha -1}\right) \left(-2 \sin \left(\frac{1}{2} \pi  (\mu+\alpha )\right)\right) & \text{ if }\varphi\in(0,\frac{\pi }{2})\,, \\
\frac{1}{2}\Gamma (1-\alpha ) (2 \varphi -\pi )^{\alpha -1}  \left(-2 \sin \left(\frac{1}{2} \pi  (\mu-\alpha )\right)\right) & \text{ if }\varphi\in(\frac{\pi }{2},\pi)\,.\\
\end{cases}
\end{equation*}
So $F_{\alpha,\beta,\mu}\in L^p([0,\pi])$ if $1<p<\frac{1}{1-\alpha}$. Then, by standard harmonic analysis \cite[page 59]{Kat04}, the convergence of the Fourier series is in $L^p([0,\pi])$. This proves \eqref{eq:integral rep exp}.

In order to prove \eqref{eq:hankel}, by Fubini's Theorem, 
$$
J_{l+m'}\left(z\right)J_{l+m}\left(z\right)=\frac{2i^{-\mu}}{\pi^2}\int_{0}^{\pi}\int_{0}^{\pi/2}e^{i2z\cos\varphi\cos\theta}(-1)^l\cos
\left((2l+\mu)\varphi\right)\cos\left(\nu\theta\right)  d\theta d\varphi.
$$
Thus, arguing as in \eqref{eq:swap int sum} we conclude that:
\begin{equation}\label{eq:almost hankel}
\sum _{l=1}^{\infty } \frac{J_{l+m'}(z)J_{l+m}(z)}{(l+\beta)^\alpha}=\frac{2i^{-\mu}}{\pi^2}\int_{0}^{\pi}\int_{0}^{\pi/2}e^{i2z\cos\varphi\cos\theta}\cos(\nu\theta)F_{\alpha,\beta,\mu}(\varphi) d\theta d\varphi.
\end{equation}
By \eqref{eq:integral rep Jn} splitting the integral from $0$ to $\pi/2$ and $\pi/2$ to $\pi$ with the change of variables $\varphi'=\pi-\varphi$ in the latter case, we obtain for $z'\in\mathbb{C}$ and even $\nu$
\begin{equation*}
\int_0^{\frac{\pi }{2}} \cos (\theta  \nu) \cos (z' \cos \theta) \, d\theta=\frac{i^\nu\pi}{2}J_\nu(z')\,,
\end{equation*}
and for odd $\nu$
\begin{equation*}
\int_0^{\frac{\pi }{2}} \cos (\theta  \nu) \sin (z' \cos \theta) \, d\theta=\frac{i^{\nu-1}\pi}{2}J_\nu(z')\,.
\end{equation*}
If we take $z=r$ real, then the RHS of \ref{eq:almost hankel} must be real. By definition, $\mu=\nu$ mod 2, so using
\begin{equation*}
\int_{0}^{\pi/2}e^{iz'\cos\theta}\cos(\nu\theta) d\theta=\int_{0}^{\pi/2}\cos(z' \cos \theta)\cos(\nu\theta) d\theta+i \int_{0}^{\pi/2}\sin(z' \cos \theta)\cos(\nu\theta) d\theta
\end{equation*}
for $z'=2z\cos\varphi$, if $\mu$ is even (and so is $\nu$), $i^{-\mu}\in\bR$, so the real part will be given by the cosine integral. The same, \emph{mutatis mutandis}, for $\mu$ odd.
 \end{proof}
\begin{remark}
	In our situation, interchanging the series and the integral sign can fail, so checking that a sufficient condition holds is necessary. For instance, as we saw in Remark \ref{rem:convergence Li alpha=0}, $F_{0,\mu}$ is well-defined, in fact, after a bit of trigonometry it equals:
	\begin{equation*}
		F_{0,\mu}(\varphi)=-\frac{\cos ((\mu+1) \varphi )}{2\cos\varphi}. 
	\end{equation*}
	But it is well-know, as we will see below, that (taking $\mu=0$)
	\begin{equation*}
		\sum _{l=1}^{\infty } J_l(r)^2 =\frac{1}{2}-\frac{J_0(r){}^2}{2}.
	\end{equation*}
	But $F_{0,0}(\varphi)=-\frac 12$, and the integral representation that the proposition above would give (if the interchange of sum and integral were allowed) is false. For instance, as we will see in the next section, the integral representation would go to zero as $r\to \infty$ but the series goes to $\frac12$. Also note that for that value of $\alpha$ our argument for a sufficient condition fails as there is no $p$ such that $1<p<\frac{1}{1-0}$.
\end{remark}
\begin{remark} Let us briefly discuss the case $z\in\mathbb{C}$. For the sake of simplicity consider $m'=\beta=0$. If $z\in\mathbb{C}$ we will have:
	\begin{align*}
	\sum _{l=1}^{\infty } \frac{J_l(z)J_{l+m}(z)}{l^\alpha}&=\frac{i^{-m}}{\pi} \int_0^\pi   \, _1\tilde{F}_2\left(1;1-\frac{m}{2},\frac{m}{2}+1;-\frac{(2z\cos \varphi)^2}{4}\right)F_{\alpha,m}(\varphi)d\varphi+\\
	&+\frac{i^{-m+1}}{\pi} \frac{z}{2} \int_0^\pi \cos \varphi\cdot \, _1\tilde{F}_2\left(1;\frac{3-m}{2},\frac{m+3}{2};-\frac{(2z\cos \varphi)^2}{4}\right)F_{\alpha,m}(\varphi)d\varphi.
	\end{align*}
	Indeed, we can use
	\begin{align*}
	\int_{0}^{\pi/2}e^{iz'\cos\theta}\cos(m\theta) d\theta&=\frac{1}{2} \pi\cdot   \, _1\tilde{F}_2\left(1;1-\frac{m}{2},\frac{m}{2}+1;-\frac{z'^2}{4}\right)+\\
	&+\frac{i\pi}{4} z'\cdot \, _1\tilde{F}_2\left(1;\frac{3-m}{2},\frac{m+3}{2};-\frac{z'^2}{4}\right)\,,
	\end{align*}
	where $_1\tilde{F}_2$ is the regularized hypergeometric function. We can check that if $z=r$ (so our series is real), we recover \eqref{eq:hankel}. Indeed, then for even $m$ and $r'=2r\cos \varphi$
	\begin{equation*}
	\sum _{l=1}^{\infty } \frac{J_l(r)J_{l+m}(r)}{l^\alpha}=\frac{i^{-m}}{\pi} \int_{0}^{\pi}  \, _1\tilde{F}_2\left(1;1-\frac{m}{2},\frac{m}{2}+1;-\frac{r'^2}{4}\right)F_{\alpha,m}(\varphi) d\varphi
	\end{equation*}
	and for odd $m$
	\begin{equation*}
	\sum _{l=1}^{\infty } \frac{J_l(r)J_{l+m}(r)}{l^\alpha}=\frac{i^{-m+1}}{\pi} \int_{0}^{\pi}\frac{r'}{2}  \, _1\tilde{F}_2\left(1;\frac{3-m}{2},\frac{m+3}{2};-\frac{r'^2}{4}\right)F_{\alpha,m}(\varphi) d\varphi.
	\end{equation*}
	Now, as the hypergeometric functions (from the definition as a series) satisfy
	$$
	\, _1\tilde{F}_2(a;b,a;x)=\, _1\tilde{F}_2(a;a,b;x)=\frac{\, _0\tilde{F}_1(;b;x)}{\Gamma (a)}
	$$
	and we have the following identity
	$$
	J_{\nu }(z)=\left(\frac{z}{2}\right)^{\nu } \, _0\tilde{F}_1\left(;\nu +1;-\frac{z^2}{4}\right),
	$$
	then, for instance,
	\begin{align*}
	\,_1\tilde{F}_2\left(1;1,1;-\frac{z^2}{4}\right)=\, _0\tilde{F}_1\left(;1;-\frac{z^2}{4}\right)&=J_0(z)\,,\\
	\left(\frac{z}{2}\right)\,_1\tilde{F}_2\left(1;2,1;-\frac{z^2}{4}\right)=\left(\frac{z}{2}\right)\, _0\tilde{F}_1\left(;2;-\frac{z^2}{4}\right)&=J_1(z)\,.
	\end{align*}
	For $z\in\mathbb{C}$, Struve functions will appear. For instance, let us consider the case $m=1$. Struve functions satisfy 
	\begin{equation*}
	\textnormal{\pmb{H}}_{\nu }(z)=\left(\frac{z}{2}\right)^{\nu +1} \, _1\tilde{F}_2\left(1;\frac{3}{2},\nu +\frac{3}{2};-\frac{z^2}{4}\right)
	\end{equation*}
	and by definition
	\begin{equation*}
	\, _1\tilde{F}_2(1;a,b;z)=\, _1\tilde{F}_2(1;b,a;z),
	\end{equation*}
	then
	\begin{align*}
	 _1\tilde{F}_2\left(1;\frac{1}{2},\frac{3}{2};-\frac{z^2}{4}\right)=\,_1\tilde{F}_2\left(1;\frac{3}{2},\frac{1}{2};-\frac{z^2}{4}\right)&=\textnormal{\pmb{H}}_{-1}(z).
	\end{align*}
	Thus,
	$$
	\sum _{l=1}^{\infty } \frac{J_l(z)J_{l+1}(z)}{l^\alpha}=\frac{1}{\pi} \int_{0}^{\pi}  \, J_1(2 z\cos\varphi) F_{\alpha,m}(\varphi) d\varphi +\frac{1}{i\pi} \int_{0}^{\pi}  \, \textnormal{\pmb{H}}_{-1}(2z\cos\varphi) F_{\alpha,m}(\varphi) d\varphi.
	$$
\end{remark}

\section{Asymptotic expansions of the Bessel series}\label{sec: asymptotics neg alpha}
Now that we have two integral representations, let us use them to calculate the asymptotic behavior of our Bessel series. Let us define:
\begin{equation*}
	\mathfrak{I}_{\alpha,\beta,\mu,\nu}(r)\coloneqq\int_0^{\pi} J_\nu(r \cos (\varphi ))F_{\alpha,\beta ,\mu}(\varphi )d\varphi.
\end{equation*}
Recall that by \eqref{eq:hankel},
\begin{equation*}
		\mathfrak{I}_{\alpha,\beta,\mu,\nu}(2r)=\pi (-1)^{m'}\sum _{l=1}^{\infty } \frac{J_{l+m'}(r)J_{l+m}(r)}{(l+\beta)^\alpha}.
\end{equation*}
\subsection{Case $\alpha\notin\bN$}
\begin{proposition}\label{prop:asymp exp alpha positive}  Let $\alpha>0$ not being an integer, then:
	\begin{equation}\label{eq:asymp exp}
	\mathfrak{I}_{\alpha,\beta,\mu,\nu}(r)=\frac{C^1_{\alpha,\mu,\nu}}{r^\alpha}+\frac{C^2_{\alpha,\beta,\mu,\nu}(r)}{r}+O_{\alpha,\beta,\mu,\nu}\left(\frac{1}{r^{\gamma}}\right)\,,
	\end{equation}	
where
\begin{align*}
	C^1_{\alpha,\mu,\nu}&\coloneqq 2 \frac{4^{\alpha -1} \Gamma (1-\alpha )\sin \left(\frac{1}{2} \pi  (\alpha +\mu)\right) \Gamma \left(\frac{\nu+\alpha }{2}\right)}{\Gamma \left(\frac{1}{2} (\nu-\alpha +2)\right)},\\
	C^2_{\alpha,\beta,\mu,\nu}(r)&\coloneqq 2\left(\cos \left(\frac{\pi  \mu }{2}\right) \zeta (\alpha ,\beta +1)+2^{-\alpha } \left(\zeta \left(\alpha ,\frac{\beta +2}{2}\right)-\zeta \left(\alpha ,\frac{\beta +1}{2}\right)\right) \sin \left(r-\frac{\pi  \nu }{2}\right)\right),
\end{align*}
where $\zeta$ is the Riemann zeta function and $\gamma\coloneqq \min\{\alpha+1, 2\}$.
\end{proposition}
\begin{proof}
	First, by Lemma \ref{lemma:F alpha m} we have that $F_{\alpha,\beta, ,\mu}(\pi-\varphi)=(-1)^\mu F_{\alpha,\beta,\mu}(\varphi)$ so using $J_\nu(-x)=(-1)^\nu J_\nu(x)$ we arrive at:
	\begin{equation*}
		\int_0^{\pi} J_\nu(r \cos (\varphi ))F_{\alpha ,\beta,\mu}(\varphi )d\varphi=2\int_0^{\pi/2} J_\nu(r \cos (\varphi ))F_{\alpha ,\beta,\mu}(\varphi )d\varphi\,,
	\end{equation*}
as $(-1)^{\nu+\mu}=1$, 	i.e., it is enough to consider the integral from 0 to $\pi/2$. The two main contributions will come from the end points of the integral so let us define $\chi$ a non-decreasing smooth function such that $\chi$ is one in a neighborhood of $\pi/2$ and vanishes  near $0$. Then,
	\begin{align}\label{eq:int Jm}
	\int_0^{\frac{\pi }{2}} J_\nu(r \cos (\varphi ))F_{\alpha ,\beta,\mu}(\varphi )&=\int_0^{c'} (1-\chi(\varphi)) J_\nu(r \cos (\varphi ))F_{\alpha ,\beta,\mu}(\varphi )  \, d\varphi  \\
	&+\int_c^{\frac{\pi }{2}} \chi(\varphi) J_\nu(r \cos (\varphi ))F_{\alpha ,\beta,\mu}(\varphi )  \, d\varphi
	\end{align}
	such that $\chi(c)=0$ and $\chi(c')=1$ with $c<c'$. Let us focus first in the second integral of the RHS. First, if $y\coloneqq \cos \varphi$, then
	$$
	\int_c^{\frac{\pi }{2}} \chi(\varphi) J_\nu(r \cos (\varphi ))F_{\alpha ,\beta,\mu}(\varphi )  \, d\varphi=\int_0^{a} \chi(\varphi(y)) J_\nu(r y)F_{\alpha ,\beta,\mu}(\varphi(y) )\frac{1}{\sqrt{1-y^2}}  \, dy \,,
	$$
	with $\varphi(y)\coloneqq\arccos y$ and  $c=\arccos(a)$. As $\chi(\varphi(a))=0$ and $a<1$, 
	\begin{equation*}
		\int_0^{a} \chi(\varphi(y)) J_\nu(r y)F_{\alpha ,\beta,\mu}(\varphi(y) )\frac{1}{\sqrt{1-y^2}}  \, dy=\int_0^{\infty} \chi(\varphi(y)) J_\nu(r y)F_{\alpha ,\beta,\mu}(\varphi(y) )\frac{1}{\sqrt{1-y^2}}  \, dy.
	\end{equation*}
	Now we are in position to analyze the integral using the asymptotic expansion of the Hankel transform \cite[Chapter IV]{Wong01}. We will need the asymptotic expansion of $F_{\alpha ,\beta,\mu}(\varphi(y))$ as $y\to0^+$. First, from the Taylor series of the exponential and the $\arccos$ function:
	\begin{equation*}
		-e^{\pm 2 i \arccos(y)}=-e^{\pm i \pi  } \pm 2 i e^{ \pm i \pi }y+2 e^{ \pm i \pi } y^2+O\left(y^3\right).
	\end{equation*}
	On the other hand, for $\alpha\neq 0,1,2,\ldots$ we can use \eqref{eq:series LerchT} with $z(y)=e^{2 i(\arccos y-\pi/2)}$. As 
	\begin{equation*}
		\arccos y-\pi/2=\sum_{n=1}^\infty a_n y^n,
	\end{equation*}
the Taylor series of $\arccos$ where $a_1=-1$. Thus,
\begin{equation*}
	\left(2i(\arccos y-\pi/2)\right)^{\alpha-1}=(-2iy)^{\alpha-1}\left(1+\sum_{n=2}^\infty (-a_n) y^{n-1}\right)^{\alpha-1}=(-2iy)^{\alpha-1}\sum_{n=0}^\infty b_n y^n \,,
\end{equation*}
using the Taylor series of $(1+x)^{\alpha-1}$.
By introducing this in \eqref{eq:series LerchT} we arrive at:
	\begin{align}\label{eq:expansion F}
		F_{\alpha ,\beta,\mu}(\varphi(y))=y^{\alpha-1 }F_{\alpha ,\beta,\mu}^1(y)+F_{\alpha ,\beta,\mu}^2(y).
	\end{align}
with $F^i_{\alpha,\beta,\mu}$ for $i=1,2$ are Taylor series given by (first term):
\begin{equation*}
	F_{\alpha ,\beta,\mu}^1(y)=2^{\alpha -1} \Gamma (1-\alpha )\sin \left(\frac{1}{2} \pi  (\alpha +\mu )\right)+O(y)
\end{equation*}
and
\begin{align*}
	F_{\alpha ,\beta,\mu}^2(y)=\cos \left(\frac{\pi  \mu }{2}\right) \zeta (\alpha ,\beta +1)+O(y).
\end{align*}
Obviously,
\begin{equation*}
	\frac{1}{\sqrt{1-y^2}}=1+\frac{y^2}{2}+O\left(y^3\right),
\end{equation*}
so if we define 
$$f(\cdot)\equiv f(\cdot,\alpha,\beta,\mu)\coloneqq \chi(\varphi(\cdot))F_{\alpha ,\beta,\mu}(\varphi(\cdot) )\frac{1}{\sqrt{1-(\cdot)^2}},$$
then we have an asymptotic expansion if $y\to0^+$ (in the sense of \cite[p.203]{Wong01}) as $\chi(\varphi(y))=1$ in a neighborhood of $y=0^+$. Furthermore, we can see that $f$ satisfies the following\footnote{Another condition is $\nu>-\alpha$, i.e., $\nu\ge 0$. Note tha this imposes no restriction as we can always interchange $m,m'$ so that $\nu$ is non-negative.}:
\begin{enumerate}[($\textnormal{H}_1$)]
	\item $f$ is smooth in $(0,\infty)$.
	\item The asymptotic expansion of $f$ is differentiable term by term.
	\item The following integrals for the $j$-th derivative are zero
	\begin{equation*}
		\int_1^{\infty} f(y)^{(j)}t^{-\frac12}e^{iry} dy=0.
	\end{equation*}
\end{enumerate}
Indeed, the first follows from Lemma \ref{lemma:F alpha m}. The second from the fact that this holds for Taylor series and so it will hold for our fractional power series, i.e.,
\begin{equation*}
	\left(\sum_{i=0}^\infty a^1_i y^{i+\alpha-1}\right)'=(\alpha-1)y^{\alpha-2}\sum_{i=0}^\infty a^1_i y^i+
	y^{\alpha-1}\sum_{i=0}^\infty a^1_i i y^{i-1}=\sum_{i=0}^\infty a^1_i (\alpha-1+i) y^{i+\alpha-2}\,,
\end{equation*}
where $	F_{\alpha ,\beta,\mu}^1(y)=\sum_{i=0}^\infty a^1_i y^{i}$. The third property follows from the fact that $\chi$ vanishes for $y>1.$ Thus, we are in position to apply \cite[Theorem 2 on p. 204, p. 207]{Wong01} to conclude that:
\begin{align}\label{eq:HankelWong pi/2}
	\int_c^{\frac{\pi }{2}} \chi(\varphi) J_\nu(r \cos (\varphi ))F_{\alpha ,\beta,\mu}(\varphi )  \, d\varphi=&\frac{4^{\alpha -1} \Gamma (1-\alpha ) \sin \left(\frac{1}{2} \pi  (\alpha +\mu )\right) \Gamma \left(\frac{\alpha +\nu }{2}\right)}{\Gamma \left(\frac{1}{2} (-\alpha +\nu +2)\right)}\frac{1}{r^\alpha}+O\left(\frac{1}{r^{\alpha+1}}\right)\nonumber\\
	&+\cos \left(\frac{\pi  \mu }{2}\right) \zeta (\alpha ,\beta +1) \frac{1}{r}+O\left(\frac{1}{r^{2}}\right).
\end{align}
On the other hand, for the first integral of \eqref{eq:int Jm} it will be better to recover the exponential integral representation \eqref{eq:integral rep exp}, i.e.,
\begin{align}
	\int_0^{c'} (1-\chi(\varphi)) J_\nu(r \cos &(\varphi ))F_{\alpha ,\beta,\mu}(\varphi)  \, d\varphi=\nonumber\\
	&=\frac{i^{-\nu}}{\pi}\int_0^{c'}\int_{0}^{\pi}e^{ir\cos \varphi\cos\theta}(1-\chi(\varphi))\cos
	\left(\nu\theta\right)F_{\alpha ,\beta,\mu}(\varphi) d\theta d\varphi.
\end{align}
 where we have used \eqref{eq:integral rep Jn}. In this way, following our first integral representation \eqref{eq:integral rep exp}, this is written as\footnote{Now $f(x,y)$ is the phase, not the amplitude as above. We do this to be in accordance with the notation of the reference we are following in each case.}
\begin{equation*}
	\int_D g(x,y)e^{i rf(x,y)} dx dy
\end{equation*}
where $D\subset\bR^2$ is a compact set, which is the ``typical'' expression for the stationary and non-stationary phase method, see \cite[(1.1) of Chapter VIII]{Wong01}. As is well-known the contribution of the asymptotic expansion will come from:
\begin{enumerate}[I)]
	\item stationary points of the phase ($\nabla f=0$),
	\item  points on
	the boundary at which a level curve of the phase is tangential to $\partial D$ and,
	\item points where $\partial D$ has some discontinuities.
\end{enumerate}
It might seem that our case is quite particular, as it satisfies I) and III) at the same time. Indeed,
$$
\nabla f=\left(-\sin (x) \cos (y),-\cos (x) \sin (y)\right)=0\Leftrightarrow x=0 \text{ and } \left(y=0 \text{ or } y=\pi\right)\,,
$$
where $x=\varphi$ and $y=\theta$ as $c'<\pi/2$, that is, the stationary points are at the corners, where the amplitude has a local (in fact, global) extremum. This situation was studied in \cite{McWo91}. Nevertheless, in our case $g(x,y)=g(x,-y)$, $g(x, \pi-y)=g(x,\pi+y)$, $g(-x,y)=g(x, y)$ and analogously for $f$. From these properties, we can use the standard theory (e.g., \cite[Chapter 7]{Hor15} or \cite[Chapter VIII]{Wong01}) of an isolated critical point in the interior of $D$. Taking this into account, we can conclude that
\begin{equation}\label{eq:approx int sta phase}
		\int_0^{c'} (1-\chi(\varphi)) J_\nu(r \cos (\varphi ))F_{\alpha ,\beta,\mu}(\varphi)  \, d\varphi=e^{-ir}\frac{c_1}{r}+e^{ir}\frac{c_2}{r}+O\left(\frac{1}{r^{2}}\right)\,,
\end{equation}
where $c_1,c_2$ depend on the values of our $g$ and $f$ (and derivatives) at the stationary points. In particular, 
\begin{align*}
	c_1\pi i^\nu&=\frac{i \pi  \cos (\pi  \nu) F_{\alpha ,\beta,\mu}(0)}{2 }=\frac{i \pi  (-1)^\nu F_{\alpha ,\beta,\mu}(0)}{2 },\\
	c_2\pi i^\nu&=\frac{-i \pi F_{\alpha ,\beta,\mu}(0)}{2 }.
\end{align*}
So, the leading term in \eqref{eq:approx int sta phase} is
\begin{equation*}
		\frac{-F_{\alpha ,\beta,\mu}(0)}{r}\left(\frac{e^{i r}+(-1)^{\nu-1}e^{-ir}}{2i^{\nu-1}}\right)=\frac{\sin \left(r-\frac{\pi  \nu}{2}\right)}{r}2^{-\alpha } \left(\zeta \left(\alpha ,\frac{\beta +2}{2}\right)-\zeta \left(\alpha ,\frac{\beta +1}{2}\right)\right)\,,		
\end{equation*}
because  
\begin{equation*}
	\Phi (-1,s,a)=\frac{\zeta \left(s,\frac{a}{2}\right)-\zeta \left(s,\frac{a+1}{2}\right)}{2^s}\,.
\end{equation*}
Thus,
\begin{align}\label{eq:stationary phase 0}
		\int_0^{c'} (1-\chi(\varphi)) J_\nu(r \cos (\varphi ))F_{\alpha,\beta,\mu}(\varphi)  \, d\varphi=&\frac{\sin \left(r-\frac{\pi  \nu}{2}\right)}{r}2^{-\alpha } \left(\zeta \left(\alpha ,\frac{\beta +2}{2}\right)-\zeta \left(\alpha ,\frac{\beta +1}{2}\right)\right)\nonumber\\ &+O\left(\frac{1}{r^2}\right)\,.
	\end{align}
The result follows from putting together \eqref{eq:HankelWong pi/2} and \eqref{eq:stationary phase 0}.
\end{proof}
\begin{remark}\label{rem:higher order terms}
	We can straightforwardly compute higher order terms in the asymptotic expansion. For the Hankel transform, we can use \cite[Theorem 2 on p. 204, p. 207]{Wong01} obtaining more terms in \eqref{eq:expansion F} using \eqref{eq:series LerchT}. For the stationary phase method we can use, for instance, \cite[Theorem 1, p. 480]{Wong01}.
\end{remark}
\begin{remark}
	Both Hankel transform and stationary phase methods are similar as the key point is the integration by parts (for instance, compare Section 2 and Section 3 of Chapter IV in \cite{Wong01}). For the stationary phase method we use the derivative of the exponential, but for the Hankel transform we use the well-known formula
	\begin{equation}
		\frac{\partial \left(z^{\nu } J_{\nu }(z)\right)}{\partial z}=z^{\nu } J_{\nu -1}(z)\,,
	\end{equation} 
	for $z\in\mathbb{C}$. Note that this formula also holds for Struve functions.
\end{remark}
We are in position to calculate the asymptotic expansion of our series:
\begin{corollary}\label{cor:series exp alpha pos} If $\mu\coloneqq m+m'$ and $\nu\coloneqq m-m'$, for $\alpha>0$ non-integer and $\beta+1>0$, we have:
	\begin{align*}\label{eq:final asymp}
		\sum _{l=1}^{\infty }& \frac{J_{l+m'}(z)J_{l+m}(z)}{(l+\beta)^\alpha}=\frac{2^{\alpha -1} \Gamma (1-\alpha )}{\Gamma \left(\frac{1}{2} (\nu-\alpha +2)\right) \Gamma \left(\frac{1}{2} (-\nu-\alpha +2)\right)}\frac{1}{r^\alpha}+O\left(\frac{1}{r^{\gamma}}\right)+\\
		&+\frac{1}{\pi r}\left(\cos \left(\frac{\pi  \nu }{2}\right) \zeta (\alpha ,\beta +1)+2^{-\alpha } \left(\zeta \left(\alpha ,\frac{\beta +2}{2}\right)-\zeta \left(\alpha ,\frac{\beta +1}{2}\right)\right) \sin \left(2r-\frac{\pi  \mu }{2}\right)\right)
	\end{align*} 
with $\gamma\coloneqq \min\{\alpha+1, 2\}$ and
where $\zeta$ is the Riemann zeta function.
\end{corollary}
\begin{remark}
	For $\beta=0$,
	\begin{equation*}
		2^{-\alpha } \left(\zeta \left(\alpha ,\frac{\beta +2}{2}\right)-\zeta \left(\alpha ,\frac{\beta +1}{2}\right)\right)=\left(2^{1-\alpha }-1\right) \zeta (\alpha ).
	\end{equation*}
\end{remark}
\begin{proof}
	The proof follows from the integral representation \eqref{eq:hankel}, i.e.,
	\begin{equation*}
		\sum _{l=1}^{\infty } \frac{J_{l+m'}(r)J_{l+m}(r)}{(l+\beta)^\alpha}=\frac{(-1)^{m'}}{\pi}\mathfrak{I}_{\alpha,\beta,\mu,\nu}(2r)
	\end{equation*}
	 and its asymptotic expansion \eqref{eq:asymp exp}, taking into account that
	$$
	\sin \left(\frac{1}{2} \pi  (\alpha +\mu)\right)
	=\sin \left(\frac{1}{2} \pi  (\alpha +\nu)\right)(-1)^{m'} \text{ , }(-1)^{{m'}} \cos \left(\frac{1}{2} \pi  \mu\right)=\cos \left(\frac{1}{2} \pi \nu\right)
	$$
	and by Euler's reflection formula we have:
	$$
	\frac{\sin \left(\pi\frac{1}{2}  (\alpha +\nu)\right) \Gamma \left(\frac{1}{2}  (\alpha +\nu)\right)}{\pi }=\frac{1}{\Gamma \left(\frac{1}{2} (-\nu-\alpha +2)\right)}.
	$$
\end{proof}
\subsection{Case $\alpha\in\bN$}
We have to analyze this case separately because \eqref{eq:series LerchT} fails. Nevertheless, from this expression we can show that:
\begin{equation}\label{eq:series LerchT integer}
\Phi(z,\alpha,v)=z^{-v} \left(\sum _{k=0,k\neq \alpha-1}^{\infty } \frac{\log ^k(z) \zeta (\alpha-k,v)}{k!}-\frac{\log ^{\alpha-1}(z) \left(\psi(v)-\psi(\alpha)+\log \left(\log \left(\frac{1}{z}\right)\right)\right)}{(\alpha-1)!}\right)
\end{equation}
where $\psi$ is the polygamma function (of order 0), see \cite[p. 30]{Bat53}. The term $\log\log (z^{-1})$ makes impossible to use the previous argument as there can be logarithmic singularities. 
\begin{proposition}\label{prop:asymp exp alpha integer} Let $\alpha\in\bN$, then:
	\begin{align*}
		\mathfrak{I}_{1,\beta,\mu,\nu}(r)=&\frac{2}{r}\left(\cos \left(\frac{\pi  \mu }{2}\right) \log (r)-\cos \left(\frac{\pi  \mu }{2}\right) \left(H_{\beta }+\psi \left(\frac{\nu +1}{2}\right)+\log (4)\right)\right)+\\
		&\frac{2}{r}\left(\frac{1}{2} \pi  \sin \left(\frac{\pi  \mu }{2}\right)-{\Phi (-1,1,\beta +1) \sin \left(r-\frac{\pi  \nu }{2}\right)}\right)	+o\left(r^{-1}\right)\,,	
	\end{align*}	
	where $H_{\beta }=\psi(\beta +1)+\gamma$, the (extended) harmonic numbers. For $\alpha>1$,
	\begin{equation*}
		\mathfrak{I}_{\alpha,\beta,\mu,\nu}(r)=\frac{C^2_{\alpha,\beta,\mu,\nu}(r)}{r}+O\left(\frac{1}{r^{2-\gamma}}\right)\,,
	\end{equation*}
where $C^2_{\alpha,\beta,\mu,\nu}(r)$ was defined in Proposition \textnormal{\ref{prop:asymp exp alpha positive}} and $\gamma>0$ as small as we want.
\end{proposition}
\begin{remark}\label{rem:higher order terms integer}
	For $\alpha=1$ we have computed second order terms of the asymptotic expansion, this will be needed in order to obtain the asymptotic expansion for $\alpha\le 0$. As in Remark \ref{rem:higher order terms}, we could do the same for higher order terms.
\end{remark}
\begin{proof}
	We proceed as in the proof of Proposition \ref{prop:asymp exp alpha positive} but now we use \eqref{eq:series LerchT integer} instead of \eqref{eq:series LerchT} and we arrive at
	\begin{align*}
		F_{\alpha ,\beta,\mu}(\varphi(y))=\log y F_{\alpha ,\beta,\mu}^1(y)+F_{\alpha ,\beta,\mu}^2(y)\,,
	\end{align*}
	with $F^i_{\alpha,\beta,\mu}$ for $i=1,2$ are Taylor series. If $z=e^\zeta$, then the second term of \eqref{eq:series LerchT integer} is
	\begin{equation*}
		\frac{\zeta ^{\alpha-1} \left(\psi(v)-\psi(\alpha)+\log \left(-\zeta\right)\right)}{(\alpha-1)!},
	\end{equation*} 
then for $\zeta$ near 0, we have that this term goes to zero if $\alpha>1$, but there is a logarithmic singularity if $\alpha=1$. Considering the latter case, proceeding as in the proof of Proposition \ref{prop:asymp exp alpha positive}, 
\begin{equation*}
	F_{\alpha=1 ,\beta,\mu}(\varphi(y))\sim\frac{1}{2} \pi  \sin \left(\frac{\pi  \mu }{2}\right)-\cos \left(\frac{\pi  \mu }{2}\right) \left(H_{\beta }+\log (2 y)\right)\,,
\end{equation*}
where $H_{\beta }=\psi(\beta +1)+\gamma$, the (extended) harmonic numbers. As we see, the dominant term is $-\cos \left(\frac{\pi  \mu }{2}\right) \log y$.\\
On the other hand, for $\alpha>1$, we have:
\begin{equation*}
		F_{\alpha ,\beta,\mu}\left(\varphi(y)\right)\sim-\frac{2^{\alpha -1} y^{\alpha -1}  \sin \left(\frac{1}{2} \pi  (\alpha -\mu )\right) \log y}{\Gamma (\alpha )}+\cos \left(\frac{\pi  \mu }{2}\right) \zeta (\alpha ,\beta +1)
\end{equation*}
so the dominant term now is $\cos \left(\frac{\pi  \mu }{2}\right) \zeta (\alpha ,\beta +1)$, as in Proposition \ref{prop:asymp exp alpha positive} for $\alpha>1$.\\

Now, we can compute the asymptotic expansion of the integrals. For that, we need to understand the behavior of the Hankel transform when logarithmic singularities are present. This is done in \cite{Wong77}. For $\alpha=1$, $F_{1,\beta,\mu}^2(y)$ can be analyzed as in the proof of Proposition \ref{prop:asymp exp alpha positive} giving rise to
\begin{equation*}
	\int_c^{\frac{\pi }{2}} \chi(\varphi) J_\nu(r \cos (\varphi ))F_{1,\beta,\mu}^2(y(\varphi))  \, d\varphi=\frac{\frac{1}{2} \pi  \sin \left(\frac{\pi  \mu }{2}\right)-\cos \left(\frac{\pi  \mu }{2}\right) \left(H_{\beta }+\log (2)\right)}{r}+O\left(r^{-2}\right).
\end{equation*}
For $F_{1,\beta,\mu}^1(y)$ we use \cite[Theorem 1, (3.12)]{Wong77} to obtain\footnote{There is a missing 2 in the denominator of (3.12) there. This is because in (5.14) it should be $2^{\alpha}$ instead of $2^{\alpha+1}$.}
\begin{align*}
\int_c^{\frac{\pi }{2}} \chi(\varphi) J_\nu(r \cos (\varphi ))\log(y(\varphi))F_{1,\beta,\mu}^1(y(\varphi))  \, d\varphi&=\frac{\cos \left(\frac{\pi  \mu }{2}\right) \log (r)}{r}\\
&-\frac{\cos \left(\frac{\pi  \mu }{2}\right) \left( \psi\left(\frac{\nu +1}{2}\right)+\log (2)\right)}{ r}+o\left(r^{-1}\right).
\end{align*}
Using the stationary phase method, as in the proof of Proposition \ref{prop:asymp exp alpha positive}, we obtain
\begin{equation*}
	\int_0^{c'} (1-\chi(\varphi)) J_\nu(r \cos (\varphi ))F_{1,\beta,\mu}(\varphi)  \, d\varphi=-\frac{\Phi (-1,1,\beta +1) \sin \left(r-\frac{\pi  \nu }{2}\right)}{r}+O(r^{-2}).
\end{equation*}
For $\alpha>1$, the first term due to the logarithmic singularity would be\footnote{This $\gamma$ is needed for $\alpha=2$. If $\alpha=3,4,\ldots$, the total error is $O\left(r^2\right)$.} $\sim \log(r)/r^\alpha=O\left(\frac{1}{r^{2-\gamma}}\right)$ $\forall\gamma >0$. Thus, as in Proposition \ref{prop:asymp exp alpha positive},
\begin{equation*}
		\int_c^{\frac{\pi }{2}} \chi(\varphi) J_\nu(r \cos (\varphi ))F_{\alpha ,\beta,\mu}(\varphi )  \, d\varphi=\cos \left(\frac{\pi  \mu }{2}\right) \zeta (\alpha ,\beta +1) \frac{1}{r}+O\left(\frac{1}{r^{2-\gamma}}\right).
\end{equation*}
where $\gamma>0$. The stationary phase method will give us the integral near $\varphi=0$ as above.	
\end{proof}
\begin{corollary}\label{cor:series expansion alpha integer} If $\mu\coloneqq m+m'$ and without loss of generality $\nu\coloneqq m-m'\ge 0$, for $\alpha>0$ integer and $\beta+1>0$, we have for $\alpha=1$,
	\begin{align*}
		\sum _{l=1}^{\infty } \frac{J_{l+m'}(z)J_{l+m}(z)}{(l+\beta)}=&\frac{1}{\pi r}\left(\cos \left(\frac{\pi  \nu }{2}\right) \log (r)-\cos \left(\frac{\pi  \nu }{2}\right) \left(H_{\beta }+\psi \left(\frac{\nu +1}{2}\right)+\log (2)\right)\right)+\\
		&\frac{1}{\pi r}\left(\frac{1}{2} \pi  \sin \left(\frac{\pi  \nu }{2}\right)-{\Phi (-1,1,\beta +1) \sin \left(2r-\frac{\pi  \mu }{2}\right)}\right)	+o\left(r^{-1}\right)
	\end{align*}
	and for $\alpha>1$
	\begin{align*}
	\sum _{l=1}^{\infty } \frac{J_{l+m'}(z)J_{l+m}(z)}{(l+\beta)^\alpha}=&\frac{1}{\pi r}\left(\cos \frac{\pi  \nu }{2} \zeta (\alpha ,\beta +1)\right)+O\left(\frac{1}{r^{2-\gamma}}\right)+\\
	&\frac{1}{\pi r}2^{-\alpha } \left(\zeta \left(\alpha ,\frac{\beta +2}{2}\right)-\zeta \left(\alpha ,\frac{\beta +1}{2}\right)\right) \sin \left(2r-\frac{\pi  \nu }{2}\right)\,,
	\end{align*} 
	where $\zeta$ is the Riemann zeta function and $\gamma>0$ arbitrary.
\end{corollary}
\begin{remark}
	As before, for $\beta=0$,
	\begin{equation*}
		\left.2^{-\alpha } \left(\zeta \left(\alpha ,\frac{\beta +2}{2}\right)-\zeta \left(\alpha ,\frac{\beta +1}{2}\right)\right)\right\vert_{\beta=0}=\left(2^{1-\alpha }-1\right) \zeta (\alpha ).
	\end{equation*}
Also,
\begin{equation*}
	\left.\Phi (-1,1,\beta +1)\right\vert_{\beta=0}=\log 2.
\end{equation*}
\end{remark}
\begin{proof}
	The proof follows the first part of the proof of Corollary \ref{cor:series exp alpha pos}.
\end{proof}
\section{Asymptotic expansions and integral representations for $\alpha$ non-negative}\label{sec: int rep asymp pos alpha}
The following lemma will be the main tool to get the integral representations and the asymptotic expansions, as it allows us to reduce the case of $\alpha\ge 0$ to the one of $\alpha<0$.
\begin{lemma}\label{lemma:trick} For $\alpha\in\bR$ and $n\in\bZ_{\ge 0}$,
	\begin{align}\label{eq:sum Salpha}
		S_{\alpha,\beta,m,m'}(r)=&\frac{r}{2}\sum_{i=0}^n(\beta-m)^i\left(S_{\alpha-1-i,\beta,m-1,m'}(r)+S_{\alpha-1-i,\beta,m+1,m'}(r)\right)+\nonumber\\
		&+(\beta-m)^{n+1} S_{\alpha-1-n,\beta,m,m'}(r)\,,
	\end{align}
where $S_{\alpha,\beta,m,m'}$	was defined in \eqref{BesselSeries}. 
\end{lemma}
\begin{proof}
	By the recurrence relations of Bessel functions,
	\begin{equation*}
		J_{l+m}(r)=\frac{r (J_{l+m-1}(r)+J_{l+m+1}(r))}{2 (l+m)},
	\end{equation*}
we get
	\begin{equation*}
J_{l+m}(r)J_{l+m'}(r)(l+\beta)^\alpha=\left(J_{l+m+1}(r)+J_{l+m-1}(r)\right)J_{l+m'}(r)\frac{r(l+\beta)^\alpha}{2(l+m)}.
	\end{equation*}
Now,
\begin{equation*}
	\frac{(l+\beta)^\alpha}{l+m}=\frac{(l+\beta)^\alpha}{l+\beta}\frac{l+\beta}{l+m}=(l+\beta)^{\alpha-1}\frac{l+\beta}{l+m}\,.
\end{equation*}
Finally, using the geometric sum,
\begin{equation*}
	\frac{l+\beta}{l+m}=\sum_{i=0}^n\left(\frac{\beta-m}{l+\beta}\right)^i+\frac{(\beta -m)^{n+1}}{(l+m)(\beta +l)^{n} }.
\end{equation*}
Thus, this gives
\begin{equation*}
	S_{\alpha,\beta,m,m'}(r)=\frac{r}{2}\sum_{i=0}^n(\beta-m)^i\left(S_{\alpha-i-1,\beta,m-1,m'}(r)+S_{\alpha-i-1,\beta,m+1,m'}(r)\right)+R_{\alpha,\beta,m,m',n}(r)\,,
\end{equation*}
where $R$ is the remainder, i.e.,
\begin{equation*}
	R_{\alpha,\beta,m,m',n}(r)=(\beta-m)^{n+1}\sum_{l=1}^\infty\frac{r}{2(l+m)} \left(J_{l+m+1}(r)+J_{l+m-1}(r)\right)J_{l+m'}(r){(\beta +l)^{\alpha-1-n}},
\end{equation*}
thus
\begin{equation*}
	R_{\alpha,\beta,m,m',n}(r)=(\beta-m)^{n+1} S_{\alpha-1-n,\beta,m,m'}(r).
\end{equation*}
\end{proof}

We are ready to prove the main theorem of this section.
\begin{proof}[Proof of Theorem \ref{prop:series expansion S}]Let us start with the case of $\alpha=0$. First, consider that $\mu$ is odd. Then, $\nu$ is odd too and $\nu\pm 1$ is even. Hence, we only have to consider the first line and the oscillatory term of the expression in Corollary \ref{cor:series expansion alpha integer} and \eqref{eq:sum Salpha} for $n=0$. The only term that survives after adding $S_{\alpha-1,\beta,m-1,m'}(r)+S_{\alpha-1,\beta,m+1,m'}(r)$ is
	\begin{equation*}
		-\frac{1}{\pi r}\left(\cos \left(\frac{\pi  (\nu+1) }{2}\right) \psi \left(\frac{\nu +2}{2}\right)+\cos \left(\frac{\pi  (\nu-1 )}{2}\right) \psi \left(\frac{\nu }{2}\right)\right)=\frac{1}{\pi r}\sin\left(\frac{\nu\pi}{2}\right)\frac{2}{\nu}\,,
	\end{equation*}
because $\psi(z+1)-\psi (z)={1}/{z}$. Thus, if $\mu$ is odd
\begin{equation*}
	S_{0,\beta,\mu,\nu}(r)=\frac{1}{\pi \nu}\sin\left(\frac{\nu\pi}{2}\right)+o(1)\,.
\end{equation*}
For $\mu$ even, we have to consider the second line of the expression in Corollary \ref{cor:series expansion alpha integer} and \eqref{eq:sum Salpha} for $n=0$. We have taken into account that $\nu\ge0$. Thus, 
\begin{equation*}
\sin \left(\frac{\pi  (\nu+1) }{2}\right)+\sin \left(\frac{\pi  |\nu-1| }{2}\right) 
\end{equation*}
will be non-zero only if $\nu=0$, giving $2$ as a result. Hence,
\begin{equation*}
	S_{0,\beta,\mu,\nu}(r)=\begin{cases*}
		\frac{1}{2} \text{ for }\nu=0\\
		0\text{ otherwise }
	\end{cases*}+o(1)\,.
\end{equation*}
For both cases we obtain the desired result by Euler's reflection formula, i.e.,
\begin{equation*}
	S_{0,\beta,\mu,\nu}(r)=\frac{1}{2 \Gamma \left(\frac{2-\nu}{2}\right) \Gamma \left(\frac{\nu+2}{2}\right)}+o(1)\,,
\end{equation*}
understood as the continuous extension. Thus, we have established, using Corollary \ref{cor:series exp alpha pos}, that for any $-\alpha\in(-1,0]$,
\begin{equation*}
S_{-\alpha,\beta,\mu,\nu}(r)=\frac{2^{\alpha -1} \Gamma (1-\alpha )}{\Gamma \left(\frac{1}{2} (\nu-\alpha +2)\right) \Gamma \left(\frac{1}{2} (-\nu-\alpha +2)\right)}\frac{1}{r^\alpha}+o\left(\frac{1}{r^{\alpha}}\right).
\end{equation*}
Therefore, for $\alpha\in(0,1]$, using Lemma \ref{eq:sum Salpha} for $S_{\alpha,\beta,m,m'}(r)=\frac{r}{2}\left(S_{\alpha-1,\beta,m-1,m'}(r)+S_{\alpha-1,\beta,m+1,m'}(r)\right)+(\beta-m) S_{\alpha-1,\beta,m,m'}(r)$,
\begin{equation*}
	S_{\alpha,\beta,\mu,\nu}(r)=\frac{2^{-\alpha -1} \Gamma (\alpha +1) r^{\alpha }}{\Gamma \left(\frac{1}{2} (\nu+\alpha +2)\right) \Gamma \left(\frac{1}{2} (-\nu+\alpha +2)\right)}+o(r^\alpha)\,,
\end{equation*}
because, using that $\Gamma(z+1)=z\Gamma(z)$,
\begin{equation*}
	\frac{1}{2} r \left(\frac{2^{-\alpha } \Gamma (\alpha ) r^{\alpha -1}}{\Gamma \left(\frac{\alpha +\nu +2}{2} \right) \Gamma \left(\frac{\alpha -\nu }{2}\right)}+\frac{2^{-\alpha } \Gamma (\alpha ) r^{\alpha -1}}{\Gamma \left(\frac{\alpha +\nu }{2}\right) \Gamma \left(\frac{\alpha -\nu +2}{2} \right)}\right)=\frac{2^{-\alpha -1} \Gamma (\alpha +1) r^{\alpha }}{\Gamma \left(\frac{\alpha -\nu +2}{2} \right) \Gamma \left(\frac{\alpha +\nu +2}{2} \right)}.
\end{equation*}
Similarly and by induction, the above equation holds for any $\alpha\ge 0$.
\end{proof}
\begin{remark}
	For $n\ge1$, by Neumann’s Addition Theorem, see \cite[(10.23.3),(10.23.4)]{Olv10},
	\begin{equation*}
		\sum _{k=0}^{2 n} (-1)^k J_k(r) J_{2 n-k}(r)+2 \sum _{l=1}^{\infty } J_l(r) J_{l+2 n}(r)=0,\quad {J_{0}}^{2}\left(z\right)+2\sum_{k=1}^{\infty}{J_{k}}^{2}\left(z\right)=1,
	\end{equation*}
which agrees with the results obtained above.
\end{remark}
 If we define, for $\alpha<0$, $$E_{\alpha,\beta,m,m'}(\varphi,\theta)\coloneqq i^{-\mu}F_{\alpha,\beta,\mu}(\varphi)\cos(\nu\theta),~~H_{\alpha,\beta,m,m'}(r,\varphi)\coloneqq (-1)^{m'}J_{\nu}(2r\cos\varphi)F_{\alpha,\beta,\mu}(\varphi)\,,$$
 and for $\alpha\ge0$ by induction
\begin{align*}
	E_{\alpha,\beta,m,m'}(r,\varphi,\theta)\coloneqq&\frac{r}{2}\left(E_{\alpha-1-i,\beta,m-1,m'}(r)+E_{\alpha-1-i,\beta,m+1,m'}(r)\right)+\nonumber\\
	&+(\beta-m)E_{\alpha-1,\beta,m,m'}(r)\,,
\end{align*}
and
\begin{align*}
	H_{{\alpha,\beta,m,m'}}(r,\varphi)\coloneqq&\frac{r}{2}\left(H_{\alpha-1-i,\beta,m-1,m'}(r)+H_{\alpha-1-i,\beta,m+1,m'}(r)\right)+\nonumber\\
	&+(\beta-m)H_{\alpha-1,\beta,m,m'}(r),
\end{align*}
then it is a straightforward consequence of Proposition \ref{prop:integral rep} and \eqref{eq:sum Salpha} the following integral representations.

\begin{proposition}\label{prop:int rep alpha nonneg} Let $\alpha\ge 0$, $\beta+1>0$ and $m,m'=0,1,\ldots$, then the series \eqref{BesselSeries} can be expressed in the following integral forms:
	\begin{equation}\label{eq:integral rep exp alpha pos}
		\sum _{l=1}^{\infty } {J_{l+m'}(r)J_{l+m}(r)}(l+\beta)^\alpha=\frac{2}{\pi^2}\int_{0}^{\pi/2}\int_{0}^{\pi}e^{i2r\cos\varphi\cos\theta}E_{\alpha,\beta,m,m'}(r,\varphi,\theta) d\varphi d\theta\,,
	\end{equation}
	as a two-dimensional exponential oscillatory integral and
	\begin{equation}\label{eq:hankel alpha pos}
		\sum _{l=1}^{\infty } {J_{l+m'}(r)J_{l+m}(r)}(l+\beta)^\alpha=\frac{1}{\pi} \int_{0}^{\pi}  \, H_{{\alpha,\beta,m,m'}}(r,\varphi) d\varphi\,,
	\end{equation}	
	as a linear combination of one-dimensional Hankel transforms, for $r\in\mathbb{R}$ where, as above, $\mu\coloneqq m+m'$ and $\nu\coloneqq m-m'$. 
\end{proposition}
	
\section{Derivative Series}\label{sec: deriv series}
In this section, building on previous results, we present a series of corollaries that significantly enhance our understanding of the summation series involving Bessel functions and their derivatives. Specifically, for different \(\alpha\), the corollaries explore various series comprising Bessel functions \(J_l(r)\), their first and second derivatives \(J'_l(r)\) and \(J''_l(r)\), and, as before, their weighted sums with the term \((l+\beta)^\alpha\), where \(l\) is a positive integer and \(\beta > -1\).

\begin{corollary}For $\alpha>-1$ we have:
	\begin{align*}
	&	\sum _{l=1}^{\infty } (l+\beta)^{\alpha } J_l(r)^2= \frac{2^{-\alpha -1} \Gamma (\alpha +1) r^{\alpha }}{\Gamma \left(\frac{\alpha }{2}+1\right)^2}+o(r^{\alpha}),\\
	&	\sum _{l=1}^{\infty } (l+\beta)^{\alpha } J_l(r) J'_l(r)= o(r^{\alpha}),\\
	&	\sum _{l=1}^{\infty } (l+\beta)^{\alpha } J'_l(r)^2= \frac{\Gamma \left(\frac{\alpha +1}{2}\right) r^{\alpha }}{4 \sqrt{\pi } \Gamma \left(\frac{\alpha }{2}+2\right)}+o(r^{\alpha}),\\
	&	\sum _{l=1}^{\infty } (l+\beta)^{\alpha } J_l(r) J''_l(r)= -\frac{\Gamma \left(\frac{\alpha +1}{2}\right) r^{\alpha }}{4 \sqrt{\pi } \Gamma \left(\frac{\alpha }{2}+2\right)}+o(r^{\alpha}),\\
	&	\sum _{l=1}^{\infty } (l+\beta)^{\alpha } J'_l(r) J''_l(r)= o(r^{\alpha}),\\
	&	\sum _{l=1}^{\infty } (l+\beta)^{\alpha } J''_l(r)^2= \frac{3\ 2^{-\alpha -5} (\alpha +2) (\alpha +4) \Gamma (\alpha +1) r^{\alpha }}{\Gamma \left(\frac{\alpha }{2}+3\right)^2}+o(r^{\alpha}).
	\end{align*}
\end{corollary}
\begin{proof}
	By the recurrence relations for Bessel functions we know
\begin{align} \label{rem:recurrence der}
\begin{split}
	J'_l(r)&=\frac{1}{2} (J_{l-1}(r)-J_{l+1}(r))=J_{l-1}( r)-\frac{l J_l( r)}{ r}\,,\\
	J''_l(r)&=\frac{1}{4} (J_{l+2}(r)+J_{l-2}(r)-2 J_l(r))=\frac{l (l+1) J_l(r)-r (J_{l-1}(r)+r J_l(r))}{r^2}.
\end{split}
\end{align}
Thus, using \eqref{rem:recurrence der} it is a straightforward consequence of Corollary \ref{cor:series exp alpha pos} and Proposition \ref{prop:series expansion S}.
\end{proof}

\begin{remark} For $\alpha\notin\mathbb{N}$ we can improve the error term as $O(r^{\alpha-\delta})$ for some $\delta>0$. Furthermore, for the second and second-last cases, the leading term vanishes. This is somehow expected as they are, respectively, the derivatives of $\frac12\sum _{l=1}^{\infty } (l+\beta)^{\alpha } J_l(r)^2$ and $\frac12 \sum _{l=1}^{\infty } (l+\beta)^{\alpha } J'_l(r)^2$. Indeed, as (see \cite[10.14.4]{Olv10})
	$$
	|J_{\ell}(r)|\le \frac{r^\ell}{2^\ell\ell!},
	$$ 	
the convergence of these series (including the ones for derivatives) is (locally) uniform, so, as standard, we can differentiate term by term. 
In this case, computing higher terms of the asymptotic expansion of \eqref{eq:asymp exp} would give the expression of the leading term.
\end{remark}
In the same way,
\begin{corollary}\label{cor:deriv series alpha -1}For $\alpha=-1$ we have
\begin{align*}
	&\sum _{l=1}^{\infty } (l+\beta)^{-1} J_l(r)^2= \frac{-H_{\beta }+\log (2 r)+\gamma }{\pi  r}+o(r^{-1}),\\
	&\sum _{l=1}^{\infty } (l+\beta)^{-1} J_l(r) J'_l(r)= -\frac{\Phi (-1,1,\beta +1) \cos (2 r)}{\pi  r}+o(r^{-1}),\\
	&\sum _{l=1}^{\infty } (l+\beta)^{-1} J'_l(r)^2= \frac{-H_{\beta }+\log (r)+\gamma -1+\log (2)}{\pi  r}+o(r^{-1}),\\
	&\sum _{l=1}^{\infty } (l+\beta)^{-1} J_l(r) J''_l(r)= \frac{\psi (\beta +1)-\log (2 r)+1}{\pi  r}+o(r^{-1}),\\
	&\sum _{l=1}^{\infty } (l+\beta)^{-1} J'_l(r) J''_l(r)= \frac{\Phi (-1,1,\beta +1) \cos (2 r)}{\pi  r}+o(r^{-1}),\\
	&\sum _{l=1}^{\infty } (l+\beta)^{-1} J''_l(r)^2= \frac{-3 \psi (\beta +1)+3 \log (r)-4+\log (8)}{3 \pi  r}+o(r^{-1}),
\end{align*}
where $\gamma$ is the Euler–Mascheroni constant and $H_\beta$ the harmonic numbers, as above.
\end{corollary}
\begin{proof}This follows from  \eqref{rem:recurrence der} and
Corollary \ref{cor:series expansion alpha integer}.
\end{proof}
\begin{corollary} For $\alpha<-1$ we have:
\begin{align*}
	&\sum _{l=1}^{\infty } (l+\beta)^{\alpha } J_l(r)^2= \frac{\zeta (-\alpha ,\beta +1)-2^{\alpha } \sin (2 r) \left(\zeta \left(-\alpha ,\frac{\beta +1}{2}\right)-\zeta \left(-\alpha ,\frac{\beta +2}{2}\right)\right)}{\pi  r}+O(r^{-2+\gamma}),\\
	&\sum _{l=1}^{\infty } (l+\beta)^{\alpha } J_l(r) J'_l(r)= -\frac{2^{\alpha } \cos (2 r) \left(\zeta \left(-\alpha ,\frac{\beta +1}{2}\right)-\zeta \left(-\alpha ,\frac{\beta +2}{2}\right)\right)}{\pi  r}+O(r^{-2+\gamma}),\\
	&\sum _{l=1}^{\infty } (l+\beta)^{\alpha } J'_l(r)^2= \frac{\zeta (-\alpha ,\beta +1)+2^{\alpha } \sin (2 r) \left(\zeta \left(-\alpha ,\frac{\beta +1}{2}\right)-\zeta \left(-\alpha ,\frac{\beta +2}{2}\right)\right)}{\pi  r}+O(r^{-2+\gamma}),\\
	&\sum _{l=1}^{\infty } (l+\beta)^{\alpha } J_l(r) J''_l(r)= -\frac{2^{\alpha } \sin (2 r) \left(\zeta \left(-\alpha ,\frac{\beta +1}{2}\right)-\zeta \left(-\alpha ,\frac{\beta +2}{2}\right)\right)-\zeta (-\alpha ,\beta +1)}{\pi  r}+O(r^{-2+\gamma}),\\
	&\sum _{l=1}^{\infty } (l+\beta)^{\alpha } J'_l(r) J''_l(r)= \frac{2^{\alpha } \cos (2 r) \left(\zeta \left(-\alpha ,\frac{\beta +1}{2}\right)-\zeta \left(-\alpha ,\frac{\beta +2}{2}\right)\right)}{\pi  r}+O(r^{-2+\gamma}),\\
	&\sum _{l=1}^{\infty } (l+\beta)^{\alpha } J''_l(r)^2= \frac{\zeta (-\alpha ,\beta +1)-2^{\alpha } \sin (2 r) \left(\zeta \left(-\alpha ,\frac{\beta +1}{2}\right)-\zeta \left(-\alpha ,\frac{\beta +2}{2}\right)\right)}{\pi  r}+O(r^{-2+\gamma}),
	\end{align*}
where $\gamma>0$ can be arbitrarily small.
\end{corollary}
\begin{proof} Again, this follows from \eqref{rem:recurrence der} and
	Corollary \ref{cor:series exp alpha pos}, \ref{cor:series expansion alpha integer}.
\end{proof}
\section*{Acknowledgments}
I would like to express my sincere gratitude to my mentors, Alberto Enciso and Daniel Peralta-Salas, for their invaluable comments and guidance during the preparation of this manuscript. 

\newpage
\bibliographystyle{siam}
\bibliography{BibArtCP}
	
\end{document}